\documentclass[11pt]{amsart}
\usepackage{amsmath,amssymb,amsthm}
\usepackage[latin1]{inputenc}
\usepackage{version,tabularx,multicol}
\usepackage{graphicx,float,psfrag}
\usepackage{stmaryrd}

\theoremstyle{plain}
\newtheorem{theo}{Theorem}[section]
 
\newtheorem{cor}[theo]{Corollary}
 
\newtheorem{prop}[theo]{Proposition}
 
\newtheorem{lem}[theo]{Lemma}
 \newtheorem{lemma}[theo]{Lemma}

\theoremstyle{remark}
\newtheorem{rem}[theo]{Remark}

\newtheorem{assumption}{Assumption}[section]

 \def\beqlb{\begin{eqnarray}}\def\eeqlb{\end{eqnarray}}
 \def\beqnn{\begin{eqnarray*}}\def\eeqnn{\end{eqnarray*}}

 \def\ar{\!\!&}

 \def\qed{\hfill$\Box$\medskip}

\newcommand{\bcen}{\begin{center}}
\newcommand{\ecen}{\end{center}}
\newcommand{\bgeqn}{\begin{equation}}
\newcommand{\edeqn}{\end{equation}}

\def\l{\left}
\def\r{\right}

 \def\ar{\!\!\!&}

\begin{document}

\title[Maximal jump and local convergence of CB processes]{The maximal jump and local convergence of continuous-state branching processes}
\date{\today}

\author{Xin He}

\address{Xin He, School of Mathematical Sciences, Beijing Normal University, Beijing 100875, P.R.CHINA}

\email{hexin@bnu.edu.cn}

\author{Zenghu Li}

\address{Zenghu Li, School of Mathematical Sciences, Beijing Normal University, Beijing 100875, P.R.CHINA}

\email{lizh@bnu.edu.cn}


\begin{abstract}
We study the distribution of the maximal jump of continuous-state branching processes.
Several exact expressions and explicit asymptotics of both the local maximal jump and the global maximal jump are obtained.
We also compare the distribution of the maximal jump and the L\'{e}vy measure to get several absolute continuity results.
Then we study local convergence of continuous-state branching processes under various conditionings.
We obtain complete results under the conditioning of large maximal jump,
and partial results under two other conditionings, which are, the conditioning of large width, and, the conditioning of large total mass.
\end{abstract}



\maketitle

\section{Introduction}\label{s:intro}

Consider a critical or subcritical continuous-state branching process (CB process) $X=(X_t)$ with the natural filtration $(\mathcal{F}_t)$. Let $\mathbf{P}_x[X\in\cdot]$ be the distribution of $X$ under the assumption of $X_0=x$, and $\mathbf{E}_x$ the corresponding expectation. Denote by $H$ the extinction time of $X$, that is,
\beqlb \label{H}
H=\inf\{s>0: X_s=0\}.
\eeqlb
Following the terminology of trees, we call $H$ the \emph{height} of $X$. Then a classical local convergence result of CB processes states that, for any $\mathcal{F}_t$-measurable bounded random variable $F$, as $r\rightarrow\infty$,
\beqlb \label{t}
\mathbf{E}_x[F| H > r]\rightarrow \frac{1}{x}\mathbf{E}_x[e^{\alpha t}X_t\,F],
\eeqlb
where $\alpha$ is specified by $\mathbf{E}_x[X_t]=e^{-\alpha t}$. See Theorem 4.1 in Li \cite{L00} for the one-dimensional version and Proposition 3.1 in Lambert \cite{La01} for the version in (\ref{t}). The convergence in (\ref{t}) implies the following statement: When conditioned to have large height, the CB process $X$ restricted to a finite time interval $[0,t]$, that is, $(X_s,s\in [0,t])$, converges weakly to $(X^*_s,s\in [0,t])$, where $X^*$ is a certain continuous-state branching process with immigration (CBI process), whose distribution is determined by that of $X$. We call the conditioning in (\ref{t}) \emph{the conditioning of large height}, and we say that under this conditioning the conditioned $X$ \emph{converges locally} to $X^*$. In the setting of superprocesses, results closely related to (\ref{t}) have actually appeared much earlier and appeared in many papers, see Section 3.3 in \cite{E00} for an introduction and several references. In the setting of L\'{e}vy trees, the tree version of (\ref{t}) has been obtained by Duquesne in \cite{D08}. Although in the literature several other conditionings have also been considered for local convergence in the continuous-state setting, they all seem to be closely related to the conditioning of large height.

However in the discrete-state setting, various conditionings have been studied for local convergence. In the seminal
paper \cite{K86}, Kesten studied local convergence of Galton-Watson trees (GW trees) under the conditioning of large height.
Since then, several other conditionings have also been considered for GW trees: the conditioning of large total
progeny, and, the conditioning of large number of leaves. Recently in \cite{AD14a,AD14b},
Abraham and Delmas provided a convenient framework to study local convergence of GW trees,
then they used this framework to prove essentially all previous results and some new ones. Specifically,
they studied the conditioning of large number of individuals with out-degree in a given set,
which includes the conditioning of large total progeny, and, the conditioning of large number of leaves
as special cases. Also very recently,
He \cite{H14} studied a new conditioning for GW trees, that is, the conditioning of large maximal out-degree.

Inspired by \cite{AD14a,AD14b,D08,H14}, naturally one would want to study local convergence of L\'{e}vy trees
under various conditionings.
Or, to avoid technicalities related to L\'{e}vy trees,
we may study local convergence of CB processes first, which is also of independent interest.
This is exactly our purpose of the present paper.
We have also obtained some results on the distribution of the maximal jump of CB processes, which seem to be interesting on their own.
Now let us explain our main results carefully in the following two paragraphs.

First in Section \ref{s:mj}, we systematically study the distribution of the maximal jump of CB processes.
Note that for most results in this section we do not exclude the supercritical case.
Our method depends crucially on stochastic equations of CB processes, and we review this topic in Section \ref{ss:SDE}.
We call $\sup_{s\in(0,t]} \Delta X_s$ the \emph{local maximal jump} for $t\in(0,\infty)$,
and $\sup_{s\in(0,\infty)} \Delta X_s$ the \emph{global maximal jump}.
We first show in Theorem \ref{trun} that the distribution of time of the first jump in a Borel set is determined by the L\'{e}vy measure of $X$
and mass processes of \emph{truncated} CB processes of $X$.
Then for the local maximal jump, Theorem \ref{lmj} expresses its distribution in terms of the L\'{e}vy measure
and the solution of an ODE. Theorem \ref{lmjr} shows that the tail of the local maximal jump
and the tail of the L\'{e}vy measure are asymptotically of the same order.
For the global maximal jump, Theorem \ref{gmj} expresses its distribution in terms of the L\'{e}vy measure
and inverse branching mechanisms of the truncated CB processes.
Then Theorem \ref{gmjr} shows that (\emph{only}) in the subcritical case, the tail of the global maximal jump
and the tail of the L\'{e}vy measure are asymptotically of the same order.
Using excursion representation of CB processes, we also get all the corresponding results under
the excursion measure in Proposition \ref{tree}, which might be useful in the study of L\'{e}vy trees.
Next we compare the distribution of the maximal jump and the L\'{e}vy measure.
In Theorem \ref{ACsub} we show that in the critical or subcritical case,
the L\'{e}vy measure and the distribution of the global maximal jump restricted to $(0,\infty)$
are absolutely continuous with respect to each other.
In Theorem \ref{AClmj} we show that in all cases,
the L\'{e}vy measure and the distribution of the local maximal jump restricted to $(0,\infty)$
are absolutely continuous with respect to each other.
In the supercritical case, the situation for the distribution of the global maximal jump is more subtle,
see Theorem \ref{ACsuper}.

Then in Section \ref{s:LC},
we study local convergence of continuous-state branching processes under various conditionings.
First we apply the conditioning of large maximal jump to CB processes,
which corresponds to the conditioning of large maximal out-degree in the discrete-state setting.
Under this conditioning, Theorem \ref{mjc} shows that in the critical case
the local limit is again $X^*$, the same CBI process appeared under the conditioning of large height.
In the subcritical case, Theorem \ref{mjs} shows that
the local limit is a certain killed CBI process $X_*$, which is different from $X^*$.
Then we consider the conditioning of large width.
We call $\sup_{s\in[0,\infty)} X_s$ the \emph{width} of $X$,
again following the terminology of trees.
To the best of our knowledge, this conditioning seems to be new, in either the continuous-state or the discrete-state setting.
Under this conditioning, Proposition \ref{w} shows that in two special critical cases
the local limit is again $X^*$.
Next we consider the conditioning of large total mass.
This conditioning is classical in the discrete-state setting (the conditioning of large total progeny),
however to the best of our knowledge, it seems to be new for local convergence in the continuous-state setting.
Under this conditioning,
Proposition \ref{tmc} shows that in the critical and stable case,
the local limit is again $X^*$.
We also study a special subcritical case in Corollary \ref{tms}, which can be reduced to the critical and stable case.
Note that in Remark \ref{conj}, we give explicit conjectures
regarding the general situation of local convergence under the conditioning of large total mass.
Finally we consider the classical conditioning of large height. Under this conditioning,
Proposition \ref{h} shows that as $r\rightarrow\infty$,
\beqlb \label{d}
\mathbf{E}_x[F| H = r]\rightarrow \frac{1}{x}\mathbf{E}_x[e^{\alpha t}X_t\,F].
\eeqlb
Inspired by a proof strategy in \cite{AD14a,AD14b}, we argue that (\ref{d}) is slightly stronger than (\ref{t}):
(\ref{d}) implies (\ref{t}) immediately, but not vice verse.
We call (\ref{t}) the \emph{tail version} of the conditioning of large height, and (\ref{d}) the \emph{density version}.
Actually under the conditioning of large total mass, we prove the density version first,
then get the tail version automatically.
However we have to admit that, in the continuous-state setting generally the density version is more restrictive then the tail version,
since the quantity in the conditioning may not have a proper density at all.

To conclude this introduction, let us mention that it seems interesting to complete our results under the conditioning of large width,
and, the conditioning of large total mass. However currently we are unable to do that.
It also seems interesting to study local convergence of L\'{e}vy trees, under the various conditionings studied in this paper.
We leave this question for future investigations.

This paper is organized as follows. In Section \ref{s:preli}, we review several basic topics in the theory of CB processes.
In Section \ref{s:mj}, we study the distribution of the maximal jump of CB processes.
Finally in Section \ref{s:LC}, we study local convergence of continuous-state branching processes under various conditionings.

\section{Preliminaries}\label{s:preli}

In this section, we review several basic topics in the theory of CB processes. In particular, we prove several lemmas which will be used in later sections.

\subsection{Continuous-state branching processes}\label{ss:cb}
This section is mainly extracted from Section 3.1 in \cite{L10}. For more details and proofs, refer to Section 3.1 in \cite{L10}.

We consider throughout the present paper a CB process $X$ with the branching mechanism
\beqlb \label{BM}
\Phi(\lambda)=\alpha \lambda+\beta \lambda^2+\int_{(0,\infty)}\pi(d\theta)(e^{-\lambda \theta}-1+\lambda \theta),
\eeqlb
where $\alpha\in \mathbb{R}$, $\beta\in\mathbb{R}_+$, and, $\pi$ is a $\sigma$-finite measure on $(0,\infty)$ satisfying
$\int_{(0,\infty)}\pi(d\theta)(\theta \wedge \theta^2)<\infty$. we exclude the trivial case of $\Phi(\lambda)\equiv0$.
Following the terminology of L\'{e}vy processes,
we call $\pi$ the \emph{L\'{e}vy measure} of $X$.
We say $\pi$ is \emph{bounded} if its support is bounded.
The branching mechanism and the corresponding CB process are called \emph{subcritical} if $\alpha>0$, \emph{critical} if $\alpha=0$,
and \emph{supercritical} if $\alpha<0$. We also use \emph{(sub)critical} to mean critical or subcritical, that is, $\alpha\geq0$.
Let $\mathbf{P}_x[X\in\cdot]$ be the distribution of $X$ under the assumption of $X_0=x$, and $\mathbf{E}_x$ the corresponding expectation.
It is well-known that $\lim_{t\rightarrow\infty}X_t=0$ a.s. in the (sub)critical case.
Also for the branching mechanism $\Phi$ given in (\ref{BM}), we have $\mathbf{E}_x[X_t]=x e^{-\alpha t}$.

It is well-known that the distribution of $X$ can be specified by $\Phi$ as follows: For $\lambda\geq 0$,
\beqlb \label{L}
\mathbf{E}_x[\exp(-\lambda X_t)]=\exp(-x v_t(\lambda)),
\eeqlb
where $v_t(\lambda)$ is the unique locally bounded nonnegative solution of
\beqlb \label{v}
v_t(\lambda)=-\int_0^t \Phi (v_s(\lambda)) ds +\lambda.
\eeqlb
It is also well-known that $X$ is Feller, so we may assume that its sample paths are rcll.
Then let $\Delta X_s=X_s-X_{s-}$ for $s\in(0,\infty)$.
For any Borel set $ A\subset (0,\infty)$, denote the time of the first jump in $A$ by $\tau_A$, that is,
\beqlb \label{tau}
\tau_A=\inf\{s>0: \Delta X_s\in A\}.
\eeqlb

We will consider the following assumptions:

\begin{assumption} \label{H0}
$\beta>0$ or $\int_{(0,\infty)} \theta\pi(d\theta)>-\alpha$.
\end{assumption}

\begin{assumption} \label{H1}
$\beta>0$ or $\int_{(0,1)}\theta\pi(d\theta)=\infty$.
\end{assumption}

\begin{assumption} \label{H2}
There is some constant $\lambda' > 0$ such that $\Phi(\lambda) > 0$ for $\lambda \geq \lambda'$ and
$$\int_{\lambda'}^\infty 1/\Phi(\lambda)d\lambda <\infty.$$
\end{assumption}

Assumption \ref{H0} Holds if and only if $\Phi(\lambda)>0$ for some $\lambda\in(0,\infty)$. Then it can be shown that as $\lambda\rightarrow\infty$, $\Phi(\lambda)\rightarrow \infty$. Refer to page 188 in \cite{B96}.
In particular, a (sub)critical branching mechanism always satisfies Assumption \ref{H0} (recall that we exclude the trivial case of $\Phi(\lambda)\equiv0$).
Also clearly Assumption \ref{H1} implies Assumption \ref{H0}.
Note that Assumption \ref{H2} implies Assumption \ref{H1}, see Corollary 3.11 in \cite{L10}. Define $\overline{v}_t=\lim_{\lambda\rightarrow\infty}v_t(\lambda)$.
Then Assumption \ref{H2} holds if and only if $\overline{v}_t<\infty$ for some and hence for all $t>0$.
Recall the definition of $H$ from (\ref{H}). Then
$$
\mathbf{P}_x[H\leq t]=\exp(-x\overline{v}_t).
$$
Under Assumption \ref{H2}, the height $H$ is finite a.s. if and only if $X$ is (sub)critical.
Under Assumption \ref{H0}, we may properly define the inverse function of $\Phi$, which we denote by $\Phi^{-1}$.
For details on $\Phi^{-1}$, refer to page 188-189 in \cite{B96}.
Note in particular that
the inverse function $\Phi^{-1}$ is a function from $[0,\infty)$ to $[q,\infty)$, where $q$ is the largest solution of $\Phi(\lambda)=0$.

\subsection{Stochastic equations of CB processes}\label{ss:SDE}

This section is extracted from Section 9.5 in \cite{L10}.
Suppose that on a suitable filtered probability space $(\Omega,\mathcal{G},\mathcal{G}_t, \mathbf{P})$,
we have a standard $\mathcal{G}_t$-Brownian motion $B_t$ and an independent $\mathcal{G}_t$-Poisson point process $P_t$
on $(0,\infty)^2$ with characteristic measure $\pi(dz)dy$.
Let $N(ds,dz,dy)$  denote the Poisson random measures on $(0,\infty)^3$ associated with $P_t$,
$\tilde{N}(ds,dz,dy)$ the compensated measure of $N(ds,dz,dy)$.

Then we may regard the CB process $X$ with the branching mechanism $\Phi$ and $X_0=x$ as the solution of the stochastic equation
$$
X_t=x-\int_0^t \alpha X_s ds+\int_0^t \sqrt{2\beta X_s} d B_s+\int_{(0,t]}\int_{(0,\infty)}\int_{(0,X_{s-}]} z \tilde{N}(ds,dz,dy).
$$

We may use stochastic equations of CB processes to prove the following result on sample paths of CB processes, which is probably known, though no reference could be found.

\begin{lemma}\label{jumps}
If $\pi(0,\infty)<\infty$, then a.s. $X$ has finite many jumps over any bounded time interval.
If $\pi(0,\infty)=\infty$, then a.s. $X$ has infinite many jumps over any nonempty open time interval, before the extinction.
\end{lemma}
\proof
Note that $\int_{(0,t]}\int_{(0,\infty)}\int_{(0,X_{s-}]} N(ds,dz,dy)$ is the number of jumps of $X$ over the time interval $(0,t]$.
We then have the following chain of relations,
\beqnn
\mathbf{E}_x\int_{(0,t]}\int_{(0,\infty)}\int_{(0,X_{s-}]} N(ds,dz,dy)
\ar =\ar\mathbf{E}_x\int_{(0,t]}\int_{(0,\infty)}\int_{(0,\infty)}\mathbf{1}\{y\leq X_{s-}\} N(ds,dz,dy)\cr
\ar =\ar\mathbf{E}_x\int_{(0,t]}\int_{(0,\infty)}\int_{(0,\infty)}\mathbf{1}\{y\leq X_{s-}\} ds \pi(dz)dy\cr
\ar =\ar\int_{(0,t]}\int_{(0,\infty)} \mathbf{E}_x[X_{s-}] ds \pi(dz)\cr
\ar =\ar x \pi(0,\infty)\int_{(0,t]} e^{-\alpha s}ds, \cr
\eeqnn
where the second identity follows from Theorem 25.22 in \cite{K02} and the last identity follows from $\mathbf{E}_x[X_{s-}]=\mathbf{E}_x[X_s]=x e^{-\alpha s}$.
So a.s. $\int_{(0,n]}\int_{(0,\infty)}\int_{(0,X_{s-}]} N(ds,dz,dy)$ is finite for all $n$ and we are done with the first statement.

For the second statement, we first argue that
for any $t_1,\,t_2,\,a$ such that $0<t_1<t_2<\infty$ and $a>0$, we have
$\int_{(t_1,t_2)}\int_{(0,\infty)}\int_{(0,a]} N(ds,dz,dy)=\infty$ a.s.
To prove this, we use the standard method based on law of large numbers.
Clearly we can find a sequence of disjoint intervals $([a_i,b_i), i\geq 1)$ such that $\pi [a_i,b_i)\geq 1$.
Then law of large numbers clearly implies that a.s.
$$
\int_{(t_1,t_2)}\int_{(0,\infty)}\int_{(0,a]} N(ds,dz,dy)\geq \sum_i \int_{(t_1,t_2)}\int_{[a_i,b_i)}\int_{(0,a]} N(ds,dz,dy)=\infty.
$$

Now we know that a.s. $\int_{(t_1,t_2]}\int_{(0,\infty)}\int_{(0,a]} N(ds,dz,dy)=\infty$ for all rational numbers $t_1,\,t_2,\,a$ such that $0<t_1<t_2<\infty$ and $a>0$. Finally note that for the CB process $X$ over any nonempty open interval $L$ before the extinction,
by right continuity of sample paths we can choose rational numbers $t_1,\,t_2,\,a$ such that $0<t_1<t_2<\infty$, $(t_1,t_2)\subset L$, and $X_s\geq a>0$ for $s\in (t_1,t_2)$.
\qed

If $0<\pi(0,\infty)<\infty$, clearly $\mathbf{E}_x\int_{(0,\infty)}\int_{(0,\infty)}\int_{(0,X_{s-}]} N(ds,dz,dy)<\infty$ if and only if $X$ is subcritical. So that we know in the subcritical case $X$ has finite many jumps over the time interval $(0,\infty)$. In the critical or supercritical case, we settle this problem of the number of jumps over $(0,\infty)$ in the following remark.

\begin{rem} If $0<\pi(0,\infty)<\infty$, we can use Lamperti representation of CB processes (see \cite{CLU09}) to study the number of jumps over $(0,\infty)$. In the (sub)critical case, a.s. $X$ has finite many jumps over the time interval $(0,\infty)$.
In the supercritical case, a.s. on the event $\{\lim_{t\rightarrow\infty}X_t=0\}$, $X$ has finite many jumps over the time interval $(0,\infty)$,
and a.s. on the event $\{\lim_{t\rightarrow\infty}X_t=\infty\}$, $X$ has infinite many jumps over the time interval $(0,\infty)$.
We do not use this result in the present paper, so we only sketch the proof here. For the L\'{e}vy process $Y$ specified by $\mathbf{E}\exp(-\lambda (Y_t-Y_0))=\exp(t\Phi(\lambda))$ and
$Y_0=x$, we stop $Y$ when it hits 0 in finite time and still denote the stopped version by $Y$. Obviously $Y$ has finite many jumps over $(0,\infty)$ if it hits 0 in finite time, otherwise it has infinite many jumps. In the (sub)critical case, $Y$ hits 0 in finite time a.s.
Then by Lamperti representation, it is not hard to argue that in the (sub)critical case $X$ has finite many jumps over $(0,\infty)$.
In the supercritical case, we may again argue by Lamperti representation. Note in particular that $Y$ does not hit 0 in finite time corresponds to $\lim_{t\rightarrow\infty}X_t=\infty$, and $Y$ hits 0 in finite time corresponds to $\lim_{t\rightarrow\infty}X_t=0$ (see e.g., Lemma 2.4 in \cite{AD12}).
\end{rem}

\subsection{Mass processes of CB processes} \label{ss:mp}

We call $(\int_0^t X_s ds,t\geq 0)$ the mass process of $X$.
By Corollary 5.17 in \cite{L10}, we have for $\lambda\geq 0$,
\beqlb \label{mp}
\mathbf{E}_x\left[ \exp\left(-\lambda\int_0^t X_s ds \right)\right]=\exp\left(-xu_t(\lambda) \right),
\eeqlb
where $u_t(\lambda)$ is the unique locally bounded nonnegative solution of
\beqlb \label{u}
u_t(\lambda)=-\int_0^t \Phi(u_s(\lambda))ds+\int_0^t \lambda ds.
\eeqlb

In general the function $u_t(\lambda)$ has no explicit expressions,
however here we give two asymptotic results on $u_t(\lambda)$, which will be useful in Section \ref{s:mj}.

\begin{lem}\label{ODElambda} For $u_t(\lambda)$, the solution of (\ref{u}), we have
$$
\frac{\partial u_t(\lambda)}{\partial \lambda}|_{\lambda=0+}=\int_0^t e^{-\alpha s} ds=\frac{1}{\alpha}(1-e^{-\alpha t}).
$$
Note that when $\alpha=0$, we agree that $\frac{1}{\alpha}(1-e^{-\alpha t})$ means $t$.
\end{lem}

\proof By considering the corresponding differential forms, we see that (\ref{u}) is equivalent to
$$
u_t(\lambda)=-\int_0^t e^{-\alpha(t-s)}\Phi_0(u_s(\lambda))ds+\int_0^t e^{-\alpha s}\lambda ds,
$$
where $\Phi_0(\lambda)=\Phi(\lambda)-\alpha \lambda$.
Since $\Phi_0$ is critical, we have
$$\frac{\partial \Phi_0(\lambda)}{\partial \lambda}|_{\lambda=0+}=0.$$
Then it is not hard to verify that
$$
\frac{\partial u_t(\lambda)}{\partial \lambda}|_{\lambda=0+}=\int_0^t e^{-\alpha s} ds=\frac{1}{\alpha}(1-e^{-\alpha t}).
$$
\qed

\begin{lem}\label{ODEt} For any $\lambda>0$, $u_t(\lambda)$ is strictly increasing with respect to $t$. If $\Phi$ satisfies Assumption \ref{H0},
then as $t\rightarrow\infty$,
$$u_t(\lambda)\rightarrow u_\infty(\lambda)=\Phi^{-1}(\lambda).$$
If $\Phi$ does not satisfy Assumption \ref{H0}, $u_t(\lambda)\rightarrow u_\infty(\lambda)=\infty$.
Finally if $\Phi$ is (sub)critical, then $u_t(0)=0=\Phi^{-1}(0)$.
\end{lem}

\proof Let $T=\inf\{t:  \Phi(u_t(\lambda))=\lambda\}$, clearly $T>0$. It is easy to see that $u_t(\lambda)$ is strictly increasing for $t\in [0,T)$.
If $T$ is finite, by the uniqueness of the solution $u_t(\lambda)$ to the equation (\ref{u}),
we see that $u_t(\lambda)=\Phi^{-1}(\lambda)$ for $t\in [T,\infty)$.
Then by (\ref{mp}) we have
$$
\mathbf{E}_x\l[\exp(-\lambda \int_0^T X_sds)\r]=\mathbf{E}_x\l[\exp(-\lambda \int_0^\infty X_sds)\r],
$$
which implies $\int_T^\infty X_sds=0$ a.s. However this is impossible since $\mathbf{E}_x[X_T]>0$ and the sample paths of $X$ are right continuous.
So $T=\infty$.

Now suppose that $\Phi$ satisfies Assumption \ref{H0} and for some $\varepsilon$ satisfying $0<\varepsilon<\Phi^{-1}(\lambda)$ and all finite $t$,
\beqlb \label{contra}
u_t(\lambda)<\Phi^{-1}(\lambda)-\varepsilon.
\eeqlb
Then by (\ref{u}), we see that for all finite $t$,
$$
\frac{\partial u_t(\lambda)}{\partial t}\geq -\Phi(\Phi^{-1}(\lambda)-\varepsilon)+\lambda>0,
$$
which implies $\lim_{t\rightarrow \infty} u_t(\lambda)=\infty$, a contradiction to (\ref{contra}).

If $\Phi$ does not satisfy Assumption \ref{H0}, it is clear that
$$
\frac{\partial u_t(\lambda)}{\partial t}\geq \lambda>0,
$$
which implies $\lim_{t\rightarrow \infty} u_t(\lambda)=\infty$.
Finally we know that
$u_t(0)=0$ by the uniqueness of the solution to (\ref{u}).
If $\Phi$ is (sub)critical, then $\Phi^{-1}(0)=0$ by the definition of $\Phi^{-1}$.
\qed

By (\ref{mp}) and monotone convergence, we see that Lemma \ref{ODEt} contains a result on the total mass $\int_0^\infty X_s ds$. For example, if $\Phi$ satisfies Assumption \ref{H0}, then for $\lambda>0$,
$$
\mathbf{E}_x\left[ \exp\left(-\lambda\int_0^\infty X_s ds \right)\right]=\exp\left(-x\Phi^{-1}(\lambda) \right).
$$
This result is also contained in Theorem VII.1 of \cite{B96}, by the well-known fact that the total mass under $\mathbf{P}_x$
and the first passage time $T(x)$ in Theorem VII.1 of \cite{B96} have the same distribution. Note that our proof here does not rely on L\'{e}vy processes.

\subsection{Excursion representation of CB processes}\label{ss:e}
This section is extracted from Section 2.4 in \cite{L12}.
Take a CB process $X$ with the branching mechanism $\Phi$, we can define an excursion measure $\mathbf{N}$ and reconstruct $X$ from excursions.
Let $\mathbb{D}_0([0,\infty),\mathbb{R}_+)$ be the subspace of $\mathbb{D}([0,\infty),\mathbb{R}_+)$, such that all paths in
$\mathbb{D}_0([0,\infty),\mathbb{R}_+)$ start from $0$ and stop upon hitting $0$. Recall the definition of $H$ from (\ref{H}).
Specifically, $\omega\in\mathbb{D}_0([0,\infty),\mathbb{R}_+)$ if and only if $\omega\in\mathbb{D}([0,\infty),\mathbb{R}_+)$,
$\omega_0=0$ and $\omega_t=0$ for $t\geq H$.
Under Assumption \ref{H1}, we may define a $\sigma$-finite measure $\mathbf{N}$ on $\mathbb{D}_0([0,\infty),\mathbb{R}_+)$ such that:

\noindent 1. $\mathbf{N}(\{\mathbf{0}\})=0$, where $\mathbf{0}$ denotes the trivial path in $\mathbb{D}([0,\infty),\mathbb{R}_+)$, that is,
$\mathbf{0}_t=0$ for any $t$.\\
2. For $t>0$, under $\mathbf{N}$ the distribution of $\omega_t$ in $(0,\infty)$ is given by the $\sigma$-finite entrance law $l_t$, see (2.2.13) in \cite{L12}. Assume that $(A^t_m,m\geq 1)$ is a partition of $(0,\infty)$ such that $l_t(A^t_m)<\infty$. If $l_t$ is finite, then let $A^t_1=(0,\infty)$ and $A^t_m=\emptyset$ for any $m\geq 2$.
\\
3. For $t>0$ and $m\geq 1$ such that $l_t(A^t_m)>0$, under the conditional probability measure $\mathbf{N}(\cdot|\omega_t\in A^t_m)$, the process $(\omega_{t+s},s\geq 0)$ is Markov with the transition kernels of the CB process $X$. \\
4. Let $N$ be a Poisson random measure on $\mathbb{D}_0([0,\infty),\mathbb{R}_+)$ with intensity $x\mathbf{N}$.
Define the process $(e_t,t\geq 0)$ by $e_0= x$ and
$$
e_t=\int_{\mathbb{D}_0([0,\infty),\mathbb{R}_+)}\omega_t N(d\omega), \quad t>0.
$$
Then $e$ is a CB process with the branching mechanism $\Phi$.

In this reconstruction of CB processes, we also notice the following fact.

\begin{lemma}\label{not} A.s. excursions in $N$ never jump at the same time.
\end{lemma}
\proof For any $t>0$, let $N^t$ be a Poisson random measure on $\mathbb{D}_0([0,\infty),\mathbb{R}_+)$ with intensity $x\mathbf{N}|_{\{\omega_t>0\}}$,
where $\mathbf{N}|_{\{\omega_t>0\}}$ is the restriction of the excursion measure $\mathbf{N}$ to the set $\{\omega_t>0\}$.
Recall the definition of Poisson random measures.
Clearly we may construct $N^t$ by considering the following partition of $\mathbb{D}_0([0,\infty),\mathbb{R}_+)$: $(\{\omega_t \in A^t_m\}, m\geq 1)$.

Now we show that a.s. excursions in $N^t$ never jump at the same time after time $t$.
Clearly it suffices to show that for any $m_1$ and $m_2$ such that $l_t(A^t_{m_1})>0$ and $l_t(A^t_{m_2})>0$, two independent excursions under the conditional probability measures $\mathbf{N}(\cdot|\omega_t\in A^t_{m_1})$ and $\mathbf{N}(\cdot|\omega_t\in A^t_{m_2})$ respectively,
never jump at the same time after time $t$.
Here $m_1$ and $m_2$ may or may not be the same.
This can be done in the usual way, by noting that the transition semigroup of the CB process $X$ is Feller.
For details, see e.g., the Remark on page 92 of \cite{RY99} and Proposition XII.1.5 in \cite{RY99}.

For an excursion to jump after time $t$, it has to be positive at time $t$.
So we see that a.s. for all $n$, excursions in $N$ never jump at the same time after time $1/n$, which means excursions in $N$ never jump at the same time.
\qed

\section{The maximal jump}\label{s:mj}

In this section we systematically study the distribution of the maximal jump of CB processes.
Consider a CB process $X$ with the branching mechanism $\Phi$ given in (\ref{BM}). Throughout this section we assume that $\pi\neq 0$.
Recall that we call $\sup_{s\in(0,t]} \Delta X_s$ the local maximal jump for $ t \in (0,\infty)$, and $\sup_{s\in(0,\infty)} \Delta X_s$
the global maximal jump.
We shall write $\sup \Delta X$ for $\sup_{s\in(0,\infty)} \Delta X_s$.

\subsection{The local maximal jump}\label{ss:lmj} We begin with an identity on the distribution of $\tau_A$ given in (\ref{tau}),
expressed in terms of the L\'{e}vy measure $\pi$ and the mass process $\int_0^t X^A_s ds$,
where $X^A$ is the so called \emph{$A$-truncated process} of $X$.
Intuitively, $X^A$ equals $X$ minus all masses produced by jumps of sizes in $A$ along with the future evolution of these masses.

\begin{theo}\label{trun} Assume that the Borel set $A\in (a,\infty)$ for some $a>0$.
Then for any $ t \in (0,\infty)$,
$$
\mathbf{P}_x\left[ \tau_A > t\right]=\mathbf{E}_x\left[ \exp\left(-\pi(A)\int_0^t X^A_s ds \right)\right],
$$
where $X^A$ is a CB process with the branching mechanism
\beqlb \label{BMr}
\Phi^A(\lambda)=\l(\alpha+\int_A \theta\pi(d\theta)\r)\lambda+\beta \lambda^2+\int_{(0,\infty)\backslash A}\pi(d\theta)(e^{-\lambda \theta}-1+\lambda \theta).
\eeqlb
It is also true that
$$
\mathbf{P}_x\left[ \tau_A =\infty\right]=\mathbf{E}_x\left[ \exp\left(-\pi(A)\int_0^\infty X^A_s ds \right)\right].
$$
\end{theo}

\proof Clearly we only need to prove the identity for finite $t$.
Recall Setion \ref{ss:SDE}. So that we may regard the CB process $X$ with the branching mechanism $\Phi$ as the solution of
$$
X_t=x-\int_0^t \alpha X_s ds+\int_0^t \sqrt{2\beta X_s} d B_s+\int_{(0,t]}\int_{(0,\infty)}\int_{(0,X_{s-}]} z \tilde{N}(ds,dz,dy).
$$
We may also regard the CB process $X^A$ with the branching mechanism $\Phi^A$ as the solution of
$$
X^A_t=x-\int_0^t \alpha^A X^A_s ds+\int_0^t \sqrt{2\beta X^A_s} d B_s+\int_{(0,t]}\int_{(0,\infty)\backslash A}\int_{(0,X^A_{s-}]} z \tilde{N}(ds,dz,dy),
$$
where $\alpha^A=\alpha+\int_A \theta\pi(d\theta)$. Let $N^A$ and $N_A$ be the restrictions of $N$ to
$(0,\infty)\times A\times(0,\infty)$ and $(0,\infty)\times {(0,\infty)\backslash A} \times(0,\infty)$, respectively.
Then notice that $X^A$ and $N^A$ are independent, since $X^A$ is ``generated" by $B$ and $N_A$.
It is obvious that $X_t=X^A_t$ if $t<\tau_A$. So for $t<\tau_A$, we have
\beqlb \label{SDEt}
X_t=x - \int_0^t \alpha X_s ds+\int_0^t \sqrt{2\beta X_s} d B_s
+\int_{(0,t]}\int_{(0,\infty)}\int_{(0,X^A_{s-}]} z \tilde{N}(ds,dz,dy).
\eeqlb
If $\tau_A\leq t$, we see that $\tau_A=\min\{s>0: \Delta X_s\in A\}$ since the sample paths of $X$ are right continuous and $A\in (a,\infty)$ for some $a>0$. So that $\Delta X_{\tau_A}\in A$, $X^A_{\tau_A}<X_{\tau_A}$, and
\beqlb \label{SDEtau}
X_{\tau_A}=x - \int_0^{\tau_A} \alpha X_s ds+\int_0^{\tau_A} \sqrt{2\beta X_s} d B_s
+\int_{(0,{\tau_A}]}\int_{(0,\infty)}\int_{(0,X^A_{s-}]} z \tilde{N}(ds,dz,dy).
\eeqlb

We argue that the two events $\{\tau_A > t\}$ and $\{\int_{(0,t]}\int_A\int_{(0,X^A_{s-}]} N(ds,dz,dy)=0\}$ coincide,
up to a null set. If $\tau_A>t$, by (\ref{SDEt}) surely $\int_{(0,t]}\int_A\int_{(0,X^A_{s-}]} N(ds,dz,dy)=0$.
If $\tau_A\leq t$, by (\ref{SDEtau}) we see that $\int_{(0,\tau_A]}\int_A\int_{(0,X^A_{s-}]} z N(ds,dz,dy)>0$,
so
$$\int_{(0,t]}\int_A\int_{(0,X^A_{s-}]} N(ds,dz,dy)>0.$$

Finally we get
\beqnn
\mathbf{P}_x\l[ \tau_A > t\r]\ar=\ar\mathbf{P}_x\l[ \int_{(0,t]}\int_A\int_{(0,X^A_{s-}]} N(ds,dz,dy)=0\r]\cr
\ar=\ar\mathbf{E}_x\left[ \exp\left(-\pi(A)\int_0^t X^A_s ds \right)\right].\cr
\eeqnn
Note that in the last step we used the independence of $X^A$ and $N^A$, where $N^A$ is the restriction of $N$ to
$(0,\infty)\times A\times(0,\infty)$.
\qed

Recall from Section \ref{ss:mp} that for $\lambda\geq 0$,
$$
\mathbf{E}_x\left[ \exp\left(-\lambda\int_0^t X^A_s ds \right)\right]=\exp\left(-x u^A_t(\lambda) \right),
$$
where $u^A_t(\lambda)$ is the unique locally bounded nonnegative solution of
\beqlb \label{ur}
u^A_t(\lambda)=-\int_0^t \Phi^A(u^A_s(\lambda))ds+\int_0^t \lambda ds.
\eeqlb
So if $A\in (a,\infty)$ for some $a>0$, Theorem \ref{trun} immediately implies
\beqlb \label{A}
\mathbf{P}_x\left[ \tau_A > t\right]=\exp\left(-x u^A_t[\pi(r,\infty)]\right).
\eeqlb
We shall write this identity in two important special cases as a theorem, since from now on essentially we only need this identity in these two special cases.

\begin{theo}\label{lmj} For $t\in(0,\infty)$ and $r\in (0,\infty)$,
$$
\mathbf{P}_x\left[ \sup_{s\in(0,t]} \Delta X_s\leq r\right]=\exp\left(-x u^{(r,\infty)}_t[\pi(r,\infty)]\right),
$$
and
$$
\mathbf{P}_x\left[ \sup_{s\in(0,t]} \Delta X_s< r\right]=\exp\left(-x u^{[r,\infty)}_t[\pi(r,\infty)]\right),
$$
where $u^{(r,\infty)}_t(\lambda)$ and $u^{[r,\infty)}_t(\lambda)$ are the unique locally bounded nonnegative solutions of (\ref{ur})
for $A=(r,\infty)$ and $A=[r,\infty)$, respectively.
\end{theo}

\proof Take $A=(r,\infty)$. Clearly $\mathbf{P}_x\left[ \sup_{s\in(0,t]} \Delta X_s\leq r\right]=\mathbf{P}_x\left[ \tau_{(r,\infty)} > t\right]$.
So the first identity follows from (\ref{A}). Similarly we get the second identity. \qed

From now on, we shall write $\Phi^r$ for $\Phi^{(r,\infty)}$, $u^r_t(\lambda)$ for $u^{(r,\infty)}_t(\lambda)$, and $\alpha^r$ for $\alpha^{(r,\infty)}$, when there is no confusion.
Using Theorem \ref{lmj} combined with Lemma \ref{ODElambda}, we can get an asymptotic result on the tail of the local maximal jump.

\begin{theo}\label{lmjr} Assume that the L\'{e}vy measure $\pi$ is unbounded. Then for any $t\in (0,\infty)$, as $r\rightarrow\infty$,
$$
\mathbf{P}_x\left[ \sup_{s\in(0,t]} \Delta X_s>r\right] \sim \frac{x}{\alpha}(1-e^{-\alpha t}) \pi(r,\infty).
$$
Note that when $\alpha=0$, we agree that $\frac{1}{\alpha}(1-e^{-\alpha t})$ means $t$.
\end{theo}

\proof Recall (\ref{BMr}). For $0<R<r<\infty$, obviously $\Phi^R\geq \Phi^r\geq \Phi$.
Then by comparing the ODEs (\ref{u}) and (\ref{ur}) we see that
$$
u^R_t(\lambda)\leq u^r_t(\lambda)\leq u_t(\lambda).
$$
Since $\pi(r,\infty)\rightarrow 0$ as $r\rightarrow\infty$, and $u_t(\lambda)\rightarrow 0$ as $\lambda\rightarrow 0$,
by Lemma \ref{ODElambda} we see that as $r\rightarrow\infty$,
$$
1-\exp\left(-x u_t[\pi(r,\infty)]\right)\sim x\pi(r,\infty)\int_0^t e^{-\alpha s} ds,
$$
and similarly
$$
1-\exp\left(-x u^R_t[\pi(r,\infty)]\right)\sim x \pi(r,\infty)\int_0^t e^{-\alpha^R s} ds,
$$
where $\alpha^R=\alpha+\int_{(R,\infty)} \theta\pi(d\theta)$.
By Theorem \ref{lmj} we see that for $0<R<r<\infty$,
$$
1-\exp\left(-x u^R_t[\pi(r,\infty)]\right)\leq \mathbf{P}_x\left[ \sup_{s\in(0,t]} \Delta X_s>r\right]\leq 1-\exp\left(-x u_t[\pi(r,\infty)]\right).
$$
We are done by finally letting $R\rightarrow\infty$ and noticing that $\lim_{R\rightarrow\infty}\alpha^R=\alpha$.
\qed

\begin{rem} From Theorem \ref{lmjr}, it is easy to see that when $\alpha\leq 0$, as $r\rightarrow\infty$,
$$
\mathbf{P}_x\left[ \sup \Delta X>r\right]/\pi(r,\infty) \rightarrow \infty,
$$
and when $\alpha> 0$, as $r\rightarrow\infty$,
$$
\mathbf{P}_x\left[ \sup \Delta X>r\right]/\pi(r,\infty) \rightarrow \frac{x}{\alpha}.
$$
However we prefer to derive this convergence in the subcritical case later in Theorem \ref{gmjr}, in the hope to make it more revealing.
\end{rem}

\subsection{The global maximal jump}\label{ss:gmj}
We begin with an identity on the distribution of the global maximal jump,
expressed in terms of the inverse function $(\Phi^r)^{-1}$ and the L\'{e}vy measure $\pi$.

\begin{theo}\label{gmj} If $\Phi$ satisfies Assumption \ref{H0}, then for any $r\in (0,\infty)$ satisfying $\pi(r,\infty)>0$,
$$
\mathbf{P}_x\l[ \sup \Delta X \leq r\r] = \exp\l(-x(\Phi^r)^{-1}[\pi(r,\infty)]\r).
$$
If $\Phi$ does not satisfy Assumption \ref{H0}, $\mathbf{P}_x\l[ \sup \Delta X \leq r\r] = 0$.
Finally if $\Phi$ is (sub)critical, then the above identity is valid even if $\pi(r,\infty)=0$.
\end{theo}

\proof
Note that $\Phi$ satisfies Assumption \ref{H0} if and only if $\Phi^r$ satisfies Assumption \ref{H0}.
Then just use Theorem \ref{lmj}, Lemma \ref{ODEt}, and monotone convergence.
\qed

\begin{rem}
In the critical and stable case, our Theorem \ref{gmj} implies Lemma 1 in \cite{B11}.
Actually, it is not hard to see that Bertoin's method in \cite{B11} can also be used to prove our Theorem \ref{gmj}.
One advantage of our method here is that we can study the local maximal jump in general,
then treat the global maximal jump as a special case.
\end{rem}

In the subcritical case, Theorem \ref{gmj} easily implies the following asymptotic result on the tail of the global maximal jump.
This result may be regarded as the continuous analogue of Theorem 3.3 in \cite{H14}, which is about subcritical GW trees.

\begin{theo}\label{gmjr} Assume that $\alpha>0$ and the L\'{e}vy measure $\pi$ is unbounded. Then as $r\rightarrow\infty$,
$$
\mathbf{P}_x\left[ \sup \Delta X>r\right] \sim \frac{x}{\alpha} \pi(r,\infty).
$$
\end{theo}

\proof The proof is very similar to that of Theorem \ref{lmjr}, however we feel that it is more revealing here.
Recall that for $0<R<r<\infty$, $\Phi^R\geq \Phi^r\geq \Phi$,
so $(\Phi^R)^{-1}\leq (\Phi^r)^{-1}\leq \Phi^{-1}$.
Since $\pi(r,\infty)\rightarrow 0$ as $r\rightarrow\infty$, $\Phi^{-1}(\lambda)\rightarrow 0$ as $\lambda\rightarrow 0$,
and
$$
\frac{\partial \Phi^{-1}(\lambda)}{\partial \lambda}|_{\lambda=0+}=1/\alpha,
$$
we see that as $r\rightarrow\infty$,
$$
1-\exp\l(-x\Phi^{-1}[\pi(r,\infty)]\r)\sim \frac{x}{\alpha}\pi(r,\infty),
$$
and similarly
$$
1-\exp\l(-x(\Phi^R)^{-1}[\pi(r,\infty)]\r)\sim \frac{x}{\alpha^R}\pi(r,\infty).
$$
By Theorem \ref{gmj} we see that for $0<R<r<\infty$,
$$
1-\exp\l(-x(\Phi^R)^{-1}[\pi(r,\infty)]\r)\leq \mathbf{P}_x\left[ \sup \Delta X>r\right]\leq 1-\exp\l(-x\Phi^{-1}[\pi(r,\infty)]\r).
$$
We are done by finally letting $R\rightarrow\infty$ and noticing that $\lim_{R\rightarrow\infty}\alpha^R=\alpha$.
\qed

By excursion representation and Lemma \ref{not}, we also get all the corresponding results under the excursion measure $\mathbf{N}$,
which might be useful in the study of L\'{e}vy trees.

\begin{prop} \label{tree} Assume that Assumption \ref{H1} holds. Then for any $t\in(0,\infty)$ and $r\in(0,\infty)$,
$$
\mathbf{N}\left[ \sup_{s\in(0,t]} \Delta \omega_s> r\right] = u^r_t[ \pi(r,\infty)].
$$
For any $r\in(0,\infty)$ satisfying $\pi(r,\infty)>0$,
$$
\mathbf{N}\left[ \sup\Delta \omega > r\right] = (\Phi^r)^{-1}[\pi(r,\infty)].
$$
Assume additionally that $\pi$ is unbounded. Then for any $t<\infty$, as $r\rightarrow\infty$,
$$
\mathbf{N}\left[ \sup_{s\in(0,t]} \Delta \omega_s>r\right] \sim \frac{1}{\alpha}(1-e^{-\alpha t}) \pi(r,\infty).
$$
Finally assume additionally that $\alpha>0$. Then as $r\rightarrow\infty$,
$$
\mathbf{N}\left[ \sup \Delta \omega >r\right] \sim \frac{1}{\alpha} \pi(r,\infty).
$$
\end{prop}

\proof First note that Assumption \ref{H1} implies Assumption \ref{H0}.
Then by excursion representation of CB processes, we see that for the global maximal jump, we have
$$
\mathbf{P}_x\left[\text{at least one excursion has the global maximal jump} > r\right]=1-\exp\l(x\mathbf{N}\l[ \sup \Delta \omega>r\r]\r).
$$
We also have a similar identity for the local maximal jump.
Recall Lemma \ref{not}. We are done by comparing the above formula to Theorem \ref{lmj}, \ref{gmj}, \ref{lmjr}, and \ref{gmjr}. \qed

\subsection{Absolute continuity of the maximal jump}\label{ss:ac} 
In this subsection, we would like to compare the distribution of the maximal jump and the L\'{e}vy measure.
First notice that the domain of the L\'{e}vy measure is $(0,\infty)$,
while it is possible for the maximal jump to be $0$, even when $\pi\neq 0$.
To discuss the possible point mass of the maximal jump at $0$,
we need the following variant of Theorem \ref{trun}.

\begin{lem} \label{finite} Assume that $\pi(0,\infty)<\infty$. Then for any $ t \in (0,\infty)$,
$$
\mathbf{P}_x\left[ \sup_{s\in(0,t]} \Delta X_s=0\right]=\mathbf{E}_x\left[ \exp\left(-\pi(0,\infty)\int_0^t X^0_s ds \right)\right],
$$
where $X^0$ is a CB process with the branching mechanism
\beqlb \label{BM0}
\Phi^0(\lambda)=\l(\alpha+\int_{(0,\infty)} \theta\pi(d\theta)\r)\lambda+\beta \lambda^2.
\eeqlb
\end{lem}

\proof We follow the proof of Theorem \ref{trun}. Define $\tau_0$ by $\tau_0=\inf\{s>0: \Delta X_s>0\}$. By Lemma \ref{jumps} we have
$\tau_0=\min\{s>0: \Delta X_s>0\}$. Then it is not hard to see that we can just follow through the proof of Theorem \ref{trun}.
\qed

From this we may also get the corresponding variants of Theorem \ref{lmj} and \ref{gmj}.
Now we can characterize the point mass of the maximal jump at $0$ completely.

\begin{prop} \label{zero} If $\pi(0,\infty)=\infty$, then for any $ t \in (0,\infty)$,
$$
\mathbf{P}_x\l[ \sup_{s\in(0,t]} \Delta X_s=0\r]=\mathbf{P}_x\l[ \sup \Delta X=0\r]=0.
$$
If $\pi(0,\infty)<\infty$, then for any $ t \in (0,\infty)$,
$$\mathbf{P}_x\left[ \sup_{s\in(0,t]} \Delta X_s=0\right]\in (0,1).$$
If $\pi(0,\infty)<\infty$ and $\Phi$ satisfies Assumption \ref{H0}, then
$$
\mathbf{P}_x\left[ \sup \Delta X=0\right]=\exp\left(-x (\Phi^0)^{-1}[\pi(0,\infty)]\right)\in (0,1),
$$
where $\Phi^0$ is given in (\ref{BM0}).
If $\pi(0,\infty)<\infty$ and $\Phi$ does not satisfy Assumption \ref{H0}, then
$\mathbf{P}_x\left[ \sup \Delta X=0\right]=0$.
\end{prop}

\proof The first statement is trivial by Lemma \ref{jumps}. Let $u^0_t(\lambda)$ be the unique locally bounded nonnegative solution of
$$
u^0_t(\lambda)=-\int_0^t \Phi^0(u^0_s(\lambda))ds+\int_0^t \lambda ds.
$$
Then if $0<\pi(0,\infty)<\infty$, by the variant of Theorem \ref{lmj} along the line of Lemma \ref{finite}, and Lemma \ref{ODEt}, we have
$$\mathbf{P}_x\left[ \sup_{s\in(0,t]} \Delta X_s=0\right]=\exp(-x u^0_t[\pi(0,\infty)])\in (0,1),$$
which is the second statement.

Notice that when $\pi(0,\infty)<\infty$, $\Phi$ satisfies Assumption \ref{H0} if and only if $\Phi^0$ does.
Then the last two statements follow from the variant of Theorem \ref{gmj} along the line of Lemma \ref{finite}.
\qed

So in general the distribution of the maximal jump is not absolutely continuous with respect to the L\'{e}vy measure, even in the (sub)critical
case. Instead we should compare the L\'{e}vy measure and the distribution of the maximal jump, restricted to $(0,\infty)$.
In the (sub)critical case, since $\lim_{t\rightarrow\infty}X_t=0$ a.s. and the sample paths of $X$ are rcll,
it is easy to see that
$\sup \Delta X=\max \Delta X$ a.s. and
$\mathbf{P}_x\l[ \sup \Delta X=\infty\r]=0$. Then Proposition \ref{zero} implies
$$\mathbf{P}_x\l[ \sup \Delta X\in(0,\infty)\r]>0.$$
Denote by $\pi|_{A}$ the restriction of the measure $\pi$ to the set $A$.
Now we compare the L\'{e}vy measure and the distribution of the global maximal jump restricted to $(0,\infty)$, in the (sub)critical case.

\begin{theo} \label{ACsub} If $\alpha\geq 0$, then $\mathbf{P}_x\l[ \sup \Delta X\in\cdot\r]|_{(0,\infty)}$
and $\pi$ are absolutely continuous with respect to each other.
\end{theo}

\proof So we need to prove $\mathbf{P}_x\l[ \sup \Delta X\in\cdot\r]|_{(0,\infty)}\ll\pi$ and $\pi \ll \mathbf{P}_x\l[ \sup \Delta X\in\cdot\r]|_{(0,\infty)}$.
The first statement is trivial, by the simple fact that $\sup \Delta X=\max \Delta X$ for $\alpha\geq 0$
and stochastic equations of CB processes. Specifically if $\pi(A)=0$ for some Borel set $A$ in $(0,\infty)$,
then by stochastic equations of CB processes clearly we have a.s. $\tau_A=\infty$, so that $\mathbf{P}_x\l[ \max \Delta X\in A\r]=0$.

For the second statement, assume that for some $0<r<r'<\infty$, $\pi(r,r')>0$. Then we are going to show that $\mathbf{P}_x\l[ \sup \Delta X\in(r,r')\r]>0$.
Now if $\mathbf{P}_x\l[ \sup \Delta X\in(r,r')\r]=0$,
then by Theorem \ref{gmj} and the variant of Theorem \ref{gmj} along the line of Theorem \ref{lmj} we have that
$$
(\Phi^{(r,\infty)})^{-1}[\pi(r,\infty)]=(\Phi^{[r',\infty)})^{-1}[\pi[r',\infty)].
$$
However this is impossible, since for any $a\geq 0$,
\beqlb \label{ele}
\Phi^{(r,\infty)}(a)=\Phi^{[r',\infty)}(a)+\int_{(r,r')}(1-e^{-a\theta})\pi(d\theta)<\Phi^{[r',\infty)}(a)+\pi(r,r').
\eeqlb
So $\mathbf{P}_x\l[ \sup \Delta X\in(r,r')\r]>0$.
Notice that $\pi(0,\infty)>0$ by our underlying assumption throughout this section and $\mathbf{P}_x\l[ \sup \Delta X\in(0,\infty)\r]>0$ by Proposition \ref{zero}.
Then it is not hard to see that for any open or closed set $A$ in $(0,\infty)$,
$\pi(A)>0$ if and only if $\mathbf{P}_x\l[ \sup \Delta X\in A\r]>0$.
Next assume that $\pi(B)>0$ for some Borel set $B$ in $(0,\infty)$.
By regularity of measures (see e.g., Lemma 1.34 in \cite{K02}), we can find some closed set $B'\subset B$, such that $\pi(B')>0$.
Then we have
$$\mathbf{P}_x\l[ \sup \Delta X\in B\r]\geq \mathbf{P}_x\l[ \sup \Delta X\in B'\r]>0,$$
which implies the second statement.
\qed

Note that in Theorem \ref{ACsub} we only consider the (sub)critical case.
The supercritical case is more subtle.
We denote by $\sup \pi$ the supremum of the support of $\pi$.

\begin{theo} \label{ACsuper}
Assume that $\alpha<0$.
Then if $\Phi$ satisfies Assumption \ref{H0},
$$
\mathbf{P}_x\left[ \sup \Delta X=\sup \pi\right]=1-\exp\left(-x (\Phi)^{-1}[0]\right)\in (0,1),
$$
and if $\Phi$ does not satisfy Assumption \ref{H0}, $\mathbf{P}_x\left[ \sup \Delta X=\sup \pi\right]=1$.
Assume additionally that $\pi(0,\sup\pi)>0$ and
$\mathbf{P}_x\l[\sup \Delta X\in(0,\sup \pi)\r]>0$, then
$\mathbf{P}_x\l[ \sup \Delta X\in\cdot\r]|_{(0,\sup\pi)}$ and $\pi|_{(0,\sup\pi)}$
are absolutely continuous with respect to each other. In particular, $\pi$ is absolutely continuous
with respect to $\mathbf{P}_x\l[ \sup \Delta X\in\cdot\r]$.
\end{theo}

\proof
For the first statement, recall the proof of Theorem \ref{gmjr}, then it is easy to see that we only need to show that
$$
\lim_{r\rightarrow\infty}(\Phi^{r})^{-1}(0)=\Phi^{-1}(0).
$$
Write $a_r$ for $(\Phi^{r})^{-1}(0)$ and $a$ for $\Phi^{-1}(0)$. Note that $a_R\leq a_r\leq a$ for $0<R<r<\infty$.
Assume that for some $R>0$ and $\varepsilon>0$, $a_r<a-\varepsilon$ for any $r>R$.
By dominated convergence we get that as $r\rightarrow\infty$,
$$
0=\Phi^{r}(a_r)\leq \Phi^{r}(a-\varepsilon)\rightarrow \Phi(a-\varepsilon)<\Phi(a)=0,
$$
a contradiction. Now we are done.

In fact, it is not hard to see that the above argument can be used to show that
under Assumption \ref{H0}, the function $(\Phi^{r})^{-1}(a)$ is continuous in $(r,a)$ over the domain $(0,\infty]\times[0,\infty)$ (Note that here we agree that $\Phi^{\infty}$ means $\Phi$),
from which the first statement is obvious. Specifically,
$$
\lim_{r\rightarrow\infty}(\Phi^{r})^{-1}[\pi(r,\infty)]=\Phi^{-1}[0].
$$

The second statement can be proved as in Theorem \ref{ACsub}.
Note in particular that by the first statement, the assumption $\mathbf{P}_x\l[\sup \Delta X\in(0,\sup \pi)\r]>0$ implies that
$\Phi$ satisfies Assumption \ref{H0}, so that we can use the identity in Theorem \ref{gmj} in the supercritical case.
Specifically, for any $r\in(0,\sup\pi)$, we have
$$
\mathbf{P}_x\l[ \sup \Delta X \leq r\r] = \exp\l(-x(\Phi^{(r,\infty)})^{-1}[\pi(r,\infty)]\r)
$$
and
$$
\mathbf{P}_x\l[ \sup \Delta X < r\r] = \exp\l(-x(\Phi^{[r,\infty)})^{-1}[\pi[r,\infty)]\r).
$$
The last statement is trivial from the first two statements.
\qed

\begin{rem}
If $\alpha<0$ and $\Phi$ does not satisfy Assumption \ref{H0},
then $\mathbf{P}_x\left[ \sup \Delta X\in (0,\sup \pi)\right]=0$, actually $\mathbf{P}_x\left[ \sup \Delta X\in [0,\sup \pi)\right]=0$.
Apart from this case, for all other cases we have that
$\mathbf{P}_x\left[ \sup \Delta X\in (0,\sup \pi)\right]>0$.
Specifically, by Proposition \ref{zero} we see that the only possible case left with $\mathbf{P}_x\left[ \sup \Delta X\in (0,\sup \pi)\right]=0$
is the case of $\pi(0,\infty)<\infty$ and $\Phi$ satisfies Assumption \ref{H0}. But for this case we still have
$\mathbf{P}_x\left[ \sup \Delta X\in (0,\sup \pi)\right]>0$,
since
$$(\Phi)^{-1}[0]<(\Phi^0)^{-1}[\pi(0,\infty)],$$
which can be verified by an inequality similar to (\ref{ele}). We can also give a counterexample of Theorem \ref{ACsub} in the supercritical case.
Clearly if $\sup \pi<\infty$ and $\pi(\sup \pi)=0$, then by Theorem \ref{ACsuper} the distribution of the global maximal jump restricted to $(0,\infty)$
is not absolutely continuous with respect to
the L\'{e}vy measure.
\end{rem}

The situation for the local maximal jump is clearer.
Since the sample paths of $X$ are rcll, we see that
a.s. $\sup_{s\in(0,t]} \Delta X_s=\max_{s\in(0,t]} \Delta X_s< \infty$ for $t\in(0,\infty)$.
Then by Proposition \ref{zero}, we see that
$$\mathbf{P}_x\l[ \sup_{s\in(0,t]} \Delta X_s\in(0,\infty)\r]>0.$$

\begin{theo} \label{AClmj}
For any $t\in(0,\infty)$, $\mathbf{P}_x\l[ \sup_{s\in(0,t]} \Delta X_s\in\cdot\r]|_{(0,\infty)}$
and $\pi$
are absolutely continuous with respect to each other.
\end{theo}

\proof
We follow the proof of Theorem \ref{ACsub}.
Again we need to prove that
$$\mathbf{P}_x\l[ \sup_{s\in(0,t]} \Delta X_s\in\cdot\r]|_{(0,\infty)}\ll\pi\quad \text{and}\quad
\pi \ll \mathbf{P}_x\l[ \sup_{s\in(0,t]} \Delta X_s\in\cdot\r]|_{(0,\infty)}.$$
As in the proof of Theorem \ref{ACsub}, the first statement is trivial.

For the second statement, assume that for some $0<r<r'<\infty$, $\pi(r,r')>0$.
Then we are going to show that $\mathbf{P}_x\l[ \sup_{s\in(0,t]} \Delta X_s\in(r,r')\r]>0$.
Now if $\mathbf{P}_x\l[ \sup_{s\in(0,t]} \Delta X_s\in(r,r')\r]=0$,
by Theorem \ref{lmj} we have
$u^{(r,\infty)}_t(\pi(r,\infty))=u^{[r',\infty)}_t(\pi[r',\infty))$.
Also by Theorem \ref{lmj} trivially $u^{(r,\infty)}_s(\pi(r,\infty))\geq u^{[r',\infty)}_s(\pi[r',\infty))$ for any $s\geq 0$.
Then it is not hard to see that $u^{(r,\infty)}_t(\pi(r,\infty))=u^{[r',\infty)}_t(\pi[r',\infty))$ implies
$$
\frac{\partial u^{(r,\infty)}_t(\pi(r,\infty))}{\partial t}=\frac{\partial u^{[r',\infty)}_t(\pi[r',\infty))}{\partial t},
$$
which by (\ref{ur}) implies that
$$
\Phi^{[r',\infty)}(a)+\pi(r,r')=\Phi^{(r,\infty)}(a)=\Phi^{[r',\infty)}(a)+\int_{(r,r')}(1-e^{-a\theta})\pi(d\theta),
$$
for $a=u^{(r,\infty)}_t(\pi(r,\infty))=u^{[r',\infty)}_t(\pi[r',\infty))$.
This is impossible by the fact that
$$\int_{(r,r')}(1-e^{-a\theta})\pi(d\theta)<\pi(r,r').$$
From here on we may just follow the proof of Theorem \ref{ACsub} to finish the present proof.
\qed

In general the distribution of $\sup \Delta X$ restricted to $(0,\infty)$ is not absolutely continuous with respect to the Lebesgue measure,
however in the (sub)critical case it is so if the L\'{e}vy measure $\pi$ is.

\begin{cor}\label{mjd} Assume that $\alpha\geq 0$ and the L\'{e}vy measure $\pi$ is absolutely continuous with respect to the Lebesgue measure.
Then the distribution of $\sup \Delta X$ under $\mathbf{P}_x$ restricted to $(0,\infty)$ has the density
function $(m_x(r),r>0)$ defined by
$$
m_x(r)=x \exp(-x(\Phi^r)^{-1}[\pi(r,\infty)])\frac{\partial (\Phi^r)^{-1}[\pi(r,\infty)]}{\partial r}.
$$
\end{cor}

\proof
The first statement follows from Theorem \ref{ACsub}.
The function $(m_x(r),r>0)$ defined here is a density of the distribution of $\sup \Delta X$ restricted to $(0,\infty)$ with respect to the Lebesgue measure,
by Theorem 3.22 in \cite{F99} and our Theorem \ref{gmj}. In particular, $\partial (\Phi^r)^{-1}[\pi(r,\infty)]/\partial r$ exists a.e. with respect to the Lebesgue measure.
\qed

Recall $\tau_A$ given in (\ref{tau}). We can also express the density function of $\tau_A$ in terms of $\Phi^A$ and $u^A_t$.

\begin{cor}
Assume that the Borel set $A\in (a,\infty)$ for some $a>0$.
Then $\tau_A$ has density $(g_A(t),t>0)$ such that
$$g_A(t)=x\l[\pi(A)-\Phi^A( u^A_t[\pi(A)])\r]\exp\l(-x u^A_t[\pi(A)]\r)\leq e^{-\alpha t} x \pi(A).$$
\end{cor}

\proof By (\ref{A}), we have
$$g_A(t)=x\l[\pi(A)-\Phi^A( u^A_t[\pi(A)])\r]\exp\l(-x u^A_t[\pi(A)]\r).$$
For the upper bound, use Theorem \ref{trun} to get
$$g_A(t)=\mathbf{E}_x\left[ \exp\left(-\pi(A)\int_{0}^{t} X^A_s ds \right)\pi(A) X^A_t\right]\leq e^{-\alpha t} x \pi(A).$$
\qed

Note that $\tau_A$ may also have a point mass at $\infty$ with the probability
$$\mathbf{P}_x[\tau_A=\infty]=\mathbf{E}_x\left[ \exp\left(-\pi(A)\int_0^\infty X^A_s ds \right)\right].$$
The r.h.s. of this identity can by evaluated by Lemma \ref{ODEt}, see the paragraph after Lemma \ref{ODEt}.

\section{Local convergence} \label{s:LC}

\subsection{Conditioning on large maximal jump}\label{ss:mj}

In this subsection we consider the (sub)critical case, that is, $\alpha\geq 0$ in (\ref{BM}).
Assume throughout this subsection that $\pi\neq 0$.
Recall that we write $\sup \Delta X$ for $\sup_{s\in(0,\infty)} \Delta X_s$.

\begin{lem}\label{simple} Assume that $\alpha\geq 0$ and the L\'{e}vy measure $\pi$ is unbounded.
Then for any $\mathcal{F}_t$-measurable bounded random variable $F$,
$$
\lim_{r\rightarrow \infty}\frac{\mathbf{E}_x\l[F \mathbf{1}\{\sup_{s\in(t,\infty)} \Delta X_s>r\}\r]}
{\mathbf{P}_x\l[\sup \Delta X >r\r]}=\frac{1}{x} \mathbf{E}_x[X_t F].
$$
\end{lem}

\proof Write $n_r=(\Phi^r)^{-1}[\pi(r,\infty)]$.
By Theorem \ref{gmj},
\beqlb \label{x}
\mathbf{P}_x[\sup \Delta X >r]=1-\exp\l(-x n_r\r).
\eeqlb
Then by Markov property,
$$
\mathbf{E}_x\l[F \mathbf{1}\{\sup_{s\in(t,\infty)} \Delta X_s > r \}\r]
=\mathbf{E}_x[F \left(1-\exp\left(-X_t n_r\right)\right)].
$$
By Theorem \ref{gmj}, we know that $\mathbf{P}_x[\sup \Delta X>r]>0$ for any $r\in(0,\infty)$.
Then by (\ref{x}) we see that $n_r>0$ for any $r\in(0,\infty)$.
In the (sub)critical case, since $\lim_{t\rightarrow\infty}X_t=0$ a.s. and the sample paths of $X$ are rcll,
it is easy to see that
a.s. $\sup \Delta X=\max \Delta X$, $\mathbf{P}_x\l[ \sup \Delta X=\infty\r]=0$, and
$$\lim_{r\rightarrow \infty}\mathbf{P}_x[\sup \Delta X>r]=0.$$
Then by (\ref{x}), we see that $\lim_{r\rightarrow \infty} n_r=0$.
Also by Theorem \ref{gmj}, we know that $\mathbf{P}_x[\sup \Delta X>r]>0$ for any $r\in(0,\infty)$.
So consequently $n_r>0$ for any $r\in(0,\infty)$.
Finally by dominated convergence,  $\mathbf{E}_x[X_t]\leq x$,
and the elementary facts that $1-e^{-a}\leq a$ for $a>0$ and $1-e^{-a}\sim a$ as $a\rightarrow 0$,
we get
$$
\lim_{r\rightarrow \infty}\frac{\mathbf{E}_x[F \left(1-\exp\left(-X_t n_r\right)\right)]}
{1-\exp\left(-x n_r\right)}=\frac{1}{x} \mathbf{E}_x[X_t F].
$$
\qed

In the critical case, Lemma \ref{simple} is already enough to imply the local convergence.

\begin{theo}\label{mjc} Assume that $\alpha = 0$ and the L\'{e}vy measure $\pi$ is unboudned.
Then for any $\mathcal{F}_t$-measurable bounded random variable $F$,
as $r\rightarrow\infty$,
$$
\mathbf{E}_x[F| \sup \Delta X >r ]\rightarrow \frac{1}{x}\mathbf{E}_x[X_t\,F].
$$
\end{theo}

\proof Let $A=\{\sup \Delta X>r\}$ and $B=\{\sup_{s\in(t,\infty)} \Delta X_s>r\}$, then
plainly we have
\beqlb \label{plain}
\mathbf{1}_A=\mathbf{1}_B+(\mathbf{1}_A-\mathbf{1}_B).
\eeqlb

For the first term on the r.h.s. of (\ref{plain}), by Lemma \ref{simple},
$$
\lim_{r\rightarrow \infty}\frac{\mathbf{E}_x[F \mathbf{1}\{\sup_{s\in(t,\infty)} \Delta X_s>r\}]}
{\mathbf{P}_x[\sup \Delta X>r]}
=\frac{1}{x} \mathbf{E}_x[X_t F].
$$
For the second term on the r.h.s. of (\ref{plain}), again by Lemma \ref{simple},
$$
\lim_{r\rightarrow \infty}\frac{\mathbf{E}_x[\mathbf{1}\{\sup \Delta X>r\}-\mathbf{1}{\{\sup_{s\in(t,\infty)} \Delta X_s>r\}}]}
{\mathbf{P}_x[\sup \Delta X >r]}=1-\frac{1}{x} \mathbf{E}_x[X_t]=0.
$$
We are done by combining the above three identities.
\qed

\begin{rem} \label{upperstar}In the (sub)critical case, define a new probability $\mathbf{P}^*_x$
by
$$
\mathbf{E}^*_x[F]=\mathbf{E}_x[e^{\alpha t}X_tF]
$$
for any $\mathcal{F}_t$-measurable bounded random variable $F$.
It is well known that under $\mathbf{P}^*_x$ the process $X$
is a CBI process with branching mechanism
$\Phi$ and immigration mechanism $\Phi'-\alpha$, where $\Phi'$ is the derivative of $\Phi$. See e.g., Section 2.3 and 3.1 in \cite{L12}
for details on CBI processes.
Use $X^*$ to denote this CBI process.
Then clearly our Theorem \ref{mjc} says that in the critical case, $X$ conditioned to have large maximal jump converges locally to $X^*$.
\end{rem}

The subcritical case is more involved, due to the fact that $1-\frac{1}{x} \mathbf{E}_x[X_t]>0$.

\begin{theo}\label{mjs} Assume that $\alpha = 0$ and the L\'{e}vy measure $\pi$ is unbounded.
Then for any $\mathcal{F}_t$-measurable bounded random variable $F$ and any $\lambda\in(0,\infty)$,
as $r\rightarrow\infty$,
$$
\mathbf{E}_x[Fe^{-\lambda X_t}| \sup \Delta X >r ]\rightarrow \frac{1}{x}\mathbf{E}_x[X_t\,Fe^{-\lambda X_t}].
$$
\end{theo}

\proof
For the first term on the r.h.s. of (\ref{plain}), again by Lemma \ref{simple}, we have
$$
\lim_{r\rightarrow \infty}\frac{\mathbf{E}_x[F \mathbf{1}\{\sup_{s\in(t,\infty)} \Delta X>r\}]}
{\mathbf{P}_x[\sup \Delta X>r]}
=\frac{1}{x} \mathbf{E}_x[X_t F].
$$
To finish the present proof, clearly it suffices to show that
$$
\lim_{r\rightarrow \infty}\frac{\mathbf{E}_x[e^{-\lambda X_t}\mathbf{1}\{\sup_{s\in(0,t]} \Delta X_s>r\}]}
{\mathbf{P}_x[\sup \Delta X>r]}=0.
$$
Recall $\tau_A$ defined in (\ref{tau}). By applying strong Markov property at $\tau_{(r,\infty)}$,
we get
\beqnn
\mathbf{E}_x\l[e^{-\lambda X_t}\mathbf{1}\{\sup_{s\in(0,t]} \Delta X_s>r\}\r]
\ar=\ar\mathbf{E}_x[e^{-\lambda X_t}\mathbf{1}\{\tau_{(r,\infty)}\leq t\}]\cr
\ar\leq \ar\mathbf{E}_x[e^{-r v_{t-{\tau_{(r,\infty)}}}(\lambda)}\mathbf{1}\{\tau_{(r,\infty)}\leq t\}]\cr
\ar\leq \ar e^{-r v_t(\lambda)}\mathbf{P}_x[\sup \Delta X>r],\cr
\eeqnn
where in the first inequality we used (\ref{L}) and the fact that $X_{\tau_{(r,\infty)}}\geq r$ if $\tau_{(r,\infty)}<\infty$,
and in the second inequality we used the fact that  $v_t(\lambda)$ is decreasing with respect to $t$, which is obvious by (\ref{v}),
and $\mathbf{P}_x[\tau_{(r,\infty)}\leq t]=\mathbf{P}_x[\sup_{s\in(0,t]} \Delta X_s>r]\leq \mathbf{P}_x[\sup \Delta X>r]$.
Finally since $v_t(\lambda)>0$ for $\lambda>0$,
we are done by letting $r\rightarrow\infty$.
\qed

\begin{rem} \label{lowerstar} In the subcritical case, note that $(\mathbf{P}_*)_x$ defined by
$$
(\mathbf{E}_*)_x[F]=\mathbf{E}_x[X_tF]
$$
is only a sub-probability. However we may extend $(\mathbf{P}_*)_x$ to a probability,
still denoted by $(\mathbf{P}_*)_x$, by letting
$$(\mathbf{P}_*)_x[X_t=\infty]=1-\mathbf{E}_x[X_t]=1-e^{-\alpha t}.$$
It is well-known that under $(\mathbf{P}_*)_x$, the process $X$
is just $X^*$ in Remark \ref{upperstar} killed at an independent exponential time with parameter $\alpha$, where by killing we mean sending $X^*$ to $\infty$.
Use $X_*$ to denote this killed process. Then by considering a special $F$ in our Theorem \ref{mjs},  it says that in the subcritical case,
$X$ conditioned to have large maximal jump converges locally to $X_*$ in the sense of finite-dimensional distributions.
We may regard $\mathbb{R}_+\cup\{\infty\}$ as a metric space homeomorphic to $[0, 1]$,
then tightness is automatic on $\mathbb{D}([0,t],\mathbb{R}_+\cup\{\infty\})$. So on $\mathbb{D}([0,t],\mathbb{R}_+\cup\{\infty\})$ weak
convergence in the sense of finite-dimensional distributions is equivalent to weak
convergence.
Finally we see that in the subcritical case, $X$ conditioned to have large maximal jump converges locally to $X_*$.
\end{rem}

\begin{rem} Assume that Assumption \ref{H1} holds, $\alpha\geq 0$, and, the L\'{e}vy measure $\pi$ is unbounded.
In this case, we can prove Lemma \ref{simple} and Theorem \ref{mjc} and \ref{mjs} without using our Theorem \ref{gmj}.
Recall Section \ref{ss:e}. By excursion representation and Lemma \ref{not}, we have
$$
\mathbf{P}_x[\sup \Delta X >r]=1-\exp\l(-x n_r\r),
$$
where $n_r=\mathbf{N} \l [\sup \Delta \omega >r \r]$.
Then it is not hard to see that we can follow through the proofs of Lemma \ref{simple} and Theorem \ref{mjc} and \ref{mjs}.
The only thing left to check is that
$\mathbf{P}_x[\sup \Delta X >r]>0$ for at least large enough $r\in (0,\infty)$.
By the one-to-one correspondence of distributions of CB processes and branching mechanisms, we know that for any $r \geq 0$,
$$
\mathbf{P}_x[\sup \Delta X >r]=\mathbf{P}_x[\tau_r <\infty]>0.
$$
\end{rem}

It should be clear that all results in this subsection are still valid in the case of $\sup \pi<\infty$ and $\pi(\sup \pi)=0$.
We just need to interpret all the limits as $r\rightarrow\sup \pi$ along the subset $\{r:r<\sup \pi\}$.

\subsection{Conditioning on large width} \label{ss:w}
In this subsection we consider two cases,
the critical and bounded L\'{e}vy measure case,
and the critical and stable case.
Recall that we call $\sup_{s\in[0,\infty)} X_s$ the width of $X$, and we write $\sup X$ for $\sup_{s\in[0,\infty)} X_s$.
We begin with a result on the tail of the width.

\begin{lemma} \label{tw} Assume that $\alpha\geq0$. Then for any $r>0$,
$$
\mathbf{P}_x[\sup X >r]\leq \frac{x}{r}.
$$
For the lower bound, assume that $\alpha=0$ and $\pi$ has bounded support,
i.e., $\pi$ has support in $[0,b]$, where $0\leq b< \infty$.
Then for any $r>x$,
$$
\mathbf{P}_x[\sup X >r]\geq\frac{x}{r+b}.
$$
So that in the critical and bounded L\'{e}vy measure case, as $r\rightarrow\infty$,
$$\mathbf{P}_x [\sup X > r] \sim \frac{x}{r}.$$
\end{lemma}

\proof The upper bound is trivial for $x\geq r$, so we assume that $r>x$ in this proof.
In the (sub)critical case, we may define $X_\infty=0$ and regard $X$ as a supermartingale over the time interval $[0,\infty]$.
Let $\tau_r=\inf\{s: X_s>r\}$ for $r\in (0,\infty)$.
By optional sampling,
$\mathbf{E}_x[X_{\tau_r}]\leq x$. Since $X_{\tau_r}\geq r$ on $\{\tau_r<\infty\}=\{\sup X >r\}$, we are done with the upper bound.

For the lower bound, note that in the critical case $X$ is a martingale over the time interval $[0,\infty)$.
We first assume that $\Phi$ satisfies Assumption \ref{H2}.
Then clearly $H \wedge \tau_r<\infty$ a.s. since $H<\infty$ a.s. by Assumption \ref{H2}.
By optional sampling, for any finite $t$,
$$x=\mathbf{E}_x[X_{H \wedge \tau_r\wedge t}]=\mathbf{E}_x[X_{H \wedge \tau_r} \mathbf{1}\{H \wedge \tau_r\leq t\}]
+\mathbf{E}_x[X_t \mathbf{1}\{H \wedge \tau_r> t\}].$$
Since $0\leq X_{H \wedge \tau_r}\leq r+b$ on $\{H \wedge \tau_r\leq t\}$ and $0\leq X_t\leq r$ on $\{H \wedge \tau_r>t\}$,
by applying dominated convergence
for both terms on the r.h.s. of the above identity, we get
$$x=\mathbf{E}_x[X_{H \wedge \tau_r}]=\mathbf{E}_x[X_{\tau_r}\mathbf{1}\{\tau_r \leq H\}]+\mathbf{E}_x[X_H \mathbf{1}\{H<\tau_r\}].$$
Now the lower bound follows from $X_H=0$ and $X_{\tau_r}\leq r+b$.

Then assume that $\Phi$ does not satisfy Assumption \ref{H2}. In this case $H=\infty$ a.s.
However it is still true that $\lim_{t\rightarrow\infty}X_t=0$ a.s.
Let $\tau'_r=\inf\{s: X_s<r\}$ for $r\in (0,\infty)$. Consider $\tau'_{1/n}$ for $n>1/x$.
Clearly $\tau'_{1/n} \wedge \tau_r<\infty$ a.s. since $\tau'_{1/n}<\infty$ a.s.
By optional sampling, for any finite $t$,
$$x=\mathbf{E}_x[X_{\tau'_{1/n} \wedge \tau_r\wedge t}]=\mathbf{E}_x[X_{\tau'_{1/n} \wedge \tau_r} \mathbf{1}\{\tau'_{1/n} \wedge \tau_r\leq t\}]
+\mathbf{E}_x[X_t \mathbf{1}\{\tau'_{1/n} \wedge \tau_r> t\}].$$
As in the previous paragraph, by dominated convergence we get
$$x=\mathbf{E}_x[X_{\tau'_{1/n} \wedge \tau_r}]=\mathbf{E}_x[X_{\tau_r}\mathbf{1}\{\tau_r \leq \tau'_{1/n}\}]+\mathbf{E}_x[X_{\tau'_{1/n}} \mathbf{1}\{\tau'_{1/n}<\tau_r\}].$$
Since $X_{\tau'_{1/n}}\leq 1/n$ and $X_{\tau_r}\leq r+b$ on $\{\tau_r \leq \tau'_{1/n}\}$,
we get
$$
\mathbf{P}_x[\tau_r \leq \tau'_{1/n}]\geq \frac{x-1/n}{r+b}.
$$
Since $\lim_{n\rightarrow\infty}\tau'_{1/n}=\infty$, letting $n\rightarrow\infty$ in the above inequality gives
$$
\mathbf{P}_x[\sup X >r]=\mathbf{P}_x[\tau_r<\infty]\geq\frac{x}{r+b}.
$$
\qed

Now we can easily derive the following result under the conditioning of large width.

\begin{prop}\label{w} Assume that $\alpha=0$ and $\pi$ has bounded support,
or $\psi(\lambda)=c\lambda^\gamma$,  where $c>0$ and $\gamma\in (1,2]$.
Then for any $\mathcal{F}_t$-measurable bounded random variable $F$, as $r\rightarrow\infty$,
$$
\mathbf{E}_x[F| \sup X >r]\rightarrow \frac{1}{x}\mathbf{E}_x[X_t\,F].
$$
\end{prop}

\proof
Trivially $\l\{\sup_{s\in [t,\infty)} X_s >r\r\} \subset  \{\sup X >r\}$,
so we get
$$\mathbf{E}_x[F| \sup X >r]
\geq \frac{\mathbf{E}_x[F\mathbf{1}\{\sup_{s\in [t,\infty)} X_s >r\}]}{\mathbf{P}_x[\sup X >r]}.
$$
In the critical and bounded L\'{e}vy measure case, by Markov property of $X$, Lemma \ref{tw} and Fatou's lemma we get
$$
\liminf_{r\rightarrow\infty} \mathbf{E}_x\l[F|\sup X >r\r] \geq \mathbf{E}_x \l[X_t F\r]/x.
$$
Clearly we may assume that $0\leq F\leq 1$, then apply the above inequality to $1-F$ to get
$$
\liminf_{r\rightarrow\infty} \mathbf{E}_x\l[1-F|\sup X >r\r] \geq \mathbf{E}_x \l[X_t (1-F)\r]/x,
$$
which implies that
$$
\limsup_{r\rightarrow\infty} \mathbf{E}_x\l[F|\sup X >r\r] \leq \mathbf{E}_x \l[X_t F\r]/x,
$$
since $\mathbf{E}_x \l[X_t\r]/x=1$. Now we are done.

In the critical and stable case, by Corollary 12.9 in \cite{K14} and explicit expressions of scale functions
(see e.g., Exercise 8.2 in \cite{K14}),
we see that as $r\rightarrow\infty$,
$$\mathbf{P}_x [\sup X > r] \sim \frac{x(\gamma-1)}{r}.$$
As in the previous paragraph, this is already enough to imply the local convergence.
\qed

By Corollary 12.9 in \cite{K14}, we can translate Lemma \ref{tw} into a result about scale functions.
We skip the details.
When $\psi(\lambda)=c\lambda^2$, we can use either Lemma \ref{tw} or the scale function method to get that for any $r\geq x$,
$$\mathbf{P}_x [\sup X > r] = \frac{x}{r}.$$

\subsection{Conditioning on large total mass}\label{ss:tm}

In this subsection we first consider the critical and stable case, that is, $\Phi(\lambda)=c\lambda^\gamma$,  where $c>0$ and $\gamma\in (1,2]$.
Denote by $\sigma$ the total mass of the branching process $X$, that is,
$$\sigma=\int_0^\infty X_s ds.$$
It is well-known that the total mass $\sigma$ under $\mathbf{P}_x$ has the same distribution as the first passage time of $Y$ below $-x$,
where $Y$ is a L\'{e}vy process specified by
$$\mathbf{E}\exp(-\lambda Y_t)=\exp(t\Phi(\lambda))=\exp(tc\lambda^\gamma).$$
Then from page 316 of \cite{DS10}
we know that under $\mathbf{P}_x$ the total mass $\sigma$ has a continuous density function
$\{f_x(t),t>0\}$.
By Remark 5 in \cite{DS10}, we get the following lemma immediately, however here we give a different proof which might be
applicable to more general cases.

\begin{lem}\label{ratio} Assume that $\psi(\lambda)=c\lambda^\gamma$,  where $c>0$ and $\gamma\in (1,2]$. Then for any $x>0$ and $y\geq 0$,
and finite $t'$,
$$
\lim_{t\rightarrow \infty}\frac{f_y(t-t')}{f_{x}(t)}=\frac{y}{x}.
$$
\end{lem}

\proof For the L\'{e}vy process $Y$ specified by $\mathbf{E}\exp(-\lambda Y_t)=\exp(tc\lambda^\gamma)$,
it is well-known that $Y_t$ has a positive continuous density $p_t(x)$
for each $t>0$, see Remark 14.18 in \cite{S99}.
Moreover, by the scaling properties, $p_t(x)$ is also continuous in $t>0$.
Then from Theorem 46.4 in \cite{S99} (or Corollary VII.3 in \cite{B96}) we have for any positive $x$ and $t$,
$$f_x(t)=\frac{x}{t}p_t(x).$$
So it suffices to prove that for any positive $x$ and $y$, and finite $t'$,
\beqlb \label{ratio limit}
\lim_{t\rightarrow \infty}\frac{p_{t-t'}(y)}{p_{t}(x)}=1.
\eeqlb
In the critical and stable case, (\ref{ratio limit}) follows easily from (14.28) on page 87 of \cite{S99}.
\qed

In the critical case, Lemma \ref{ratio} is enough to imply the local convergence.

\begin{prop}\label{tmc} Assume that $\psi(\lambda)=c\lambda^\gamma$,  where $c>0$ and $\gamma\in (1,2]$.
Then for any $\mathcal{F}_t$-measurable bounded random variable $F$,
as $r\rightarrow\infty$,
$$
\mathbf{E}_x[F| \sigma = r]\rightarrow \frac{1}{x}\mathbf{E}_x[X_t\,F] \quad \text{and} \quad \mathbf{E}_x[F| \sigma > r]
\rightarrow \frac{1}{x}\mathbf{E}_x[X_t\,F].
$$
\end{prop}

\proof First of all,
$$
\mathbf{E}_x[F| \sigma > r]=\int_r^\infty \mathbf{E}_x[F| \sigma = a] f_x(a|r)da,
$$
where $f_x(a|r)$ is the conditioned density defined by
$$
f_x(a|r)=\frac{f_x(a)}{\int_r^\infty f_x(a) da}.
$$
Since  $\int_r^\infty f_x(a|r)da=1$, clearly the local convergence of the tail version follows from that of the density version.

To prove the local convergence of the density version, we denote $\int_0^t X_s ds$ by $\sigma_t$, clearly we have
$$
\l\{\int_t^\infty X_s ds = \sigma-\sigma_t = r-\sigma_t, \sigma_t<r\r\} \subset  \{\sigma = r\}.
$$
Combined with Disintegration theorem (see e.g., Theorem 6.4 in \cite{K02}) and Markov property, we have
$$
\mathbf{E}_x\l[F | \sigma = r\r]
\geq\mathbf{E}_x\l[F\mathbf{1}\{\sigma_t<r\}f_{X_t}(r-\sigma_t)\r]/f_x(r).
$$
Then by Lemma \ref{ratio} and Fatou's lemma we get
$$
\liminf_{r\rightarrow\infty} \mathbf{E}_x\l[F|\sigma = r\r] \geq \mathbf{E}_x \l[X_t F\r]/x.
$$
Finally we can just follow the end of the proof of Proposition \ref{w} to improve the above to
$$
\lim_{r\rightarrow\infty} \mathbf{E}_x\l[F|\sigma = r\r] = \mathbf{E}_x \l[X_t F\r]/x.
$$
\qed

We then consider a special subcritical case, which can be reduced to the critical and stable case.
To state the following result, we need to introduce the shifted branching mechanisms.
Consider the branching mechanism $\Phi$ given in (\ref{BM}). We use $\Theta^\Phi$ to denote all $\theta\in \mathbb{R}$
such that $\int_1^\infty a e^{-\theta a}\pi (da)<\infty$. For any $\theta\in\Theta^\Phi$,
define a function $\Phi_\theta(\lambda)$ on $\mathbb{R}_+$ by $\Phi_\theta(\lambda)=\Phi(\theta+\lambda)-\Phi(\theta)$.
It is easy to see that the function $\Phi_\theta$ is also a branching mechanism, and
$$
\Phi_\theta(\lambda)=\alpha_\theta \lambda+\beta \lambda^2+\int_{(0,\infty)}\pi_\theta(da)(e^{-\lambda a}-1+\lambda a),
$$
where $\alpha_\theta=\alpha + 2\beta\theta+\int_0^\infty(1-e^{-\theta a})a\pi(da)$ and $\pi_\theta(da)=e^{-\theta a}\pi(da)$.

\begin{cor}\label{tms} Assume that $\alpha<0$ and there exists a negative $q\in \Theta^\Phi$ such that
$$
 \Phi_q(\lambda)=c\lambda^\gamma,
$$
where $c>0$ and $\gamma\in (1,2]$.
Then for any $\mathcal{F}_t$-measurable bounded random variable $F$,
as $r\rightarrow\infty$,
$$
\mathbf{E}_x[F| \sigma = r]\rightarrow \frac{1}{x}\mathbf{E}_x[X_tF\,e^{qx-qX_t-\Phi(q)\tau_t}]
\quad \text{and} \quad
\mathbf{E}_x[F| \sigma > r]\rightarrow \frac{1}{x}\mathbf{E}_x[X_tF\,e^{qx-qX_t-\Phi(q)\tau_t}],
$$
where $\sigma_t=\int_0^t X_s ds$.
\end{cor}

\proof As in the proof of Proposition \ref{tmc}, we only need to prove the density version.
Let $\mathbf{P}^\Phi_x$ be the law of $X$ with the branching mechanism $\Phi$ and $X_0=x$,
and $\mathbf{E}^\Phi_x$ the corresponding expectation.
Then we recall the following conditional equivalence result,
$$
\mathbf{E}_x^\Phi [  \cdot | \sigma = r ]=\mathbf{E}_x^{\Phi_q} [  \cdot | \sigma = r ],
$$
which is implied by Lemma 2.4.(ii) in \cite{AD12}. Also recall from Theorem 2.2.(ii) in \cite{AD12}
that for any $\mathcal{F}_t$-measurable bounded random variable $F$,
$$
\mathbf{E}_x^{\Phi_q} [F]=\mathbf{E}_x^\Phi [F\,e^{qx-qX_t-\Phi(q)\tau_t}].
$$
We are done with the density version by combining together the above two identities and Proposition \ref{tmc}.
\qed

\begin{rem} From the proof of Corollary \ref{tms}, we see that in the setting of this corollary, the conditioned CB process $X$
converges locally to $(X^{\Phi_q})^*$, where $(X^{\Phi_q})^*$ is a CBI process
with branching mechanism
$\Phi_q$ and immigration mechanism $\Phi'_q$. Note that $\Phi_q$ is critical.
\end{rem}

\begin{rem}\label{conj} Inspired by the corresponding results of GW trees (see e.g., Definition 1.1 and Theorem 1.3 in \cite{AD14b}),
we make the following conjectures on the general situation of local convergence under the conditioning of large total mass:\\
\noindent Case I, $\alpha=0$, then the conditioned $X$ converges locally to $X^*$, where $X^*$ is a CBI process
with branching mechanism $\Phi$ and immigration mechanism $\Phi'$;\\
Case II, $\alpha<0$ and there exists a negative $q\in \Theta^\Phi$ such that
$\Phi_q$ is critical, then the conditioned $X$ converges locally to $(X^{\Phi_q})^*$, where $(X^{\Phi_q})^*$ is a CBI process
with branching mechanism $\Phi_q$ and immigration mechanism $\Phi'_q$;\\
Case III, $\alpha<0$ and $\Phi_q$ is subcritical for any $q\in \Theta^\Phi$,
then $\inf\Theta^\Phi\in \Theta^\Phi$ (easy to check) and we denote it by $q'$. Note that $\Phi_{q'}$ is also subcritical.
Finally
the conditioned $X$ converges locally to $(X^{\Phi_{q'}})_*$, where $(X^{\Phi_{q'}})_*$ is a killed CBI process
with branching mechanism $\Phi_{q'}$ and immigration mechanism $\Phi'_{q'}-\alpha_{q'}$,
and killed at an independent exponential time with parameter $\alpha_{q'}$. One may refer to Remark \ref{upperstar} and \ref{lowerstar} for more details on this killed CBI process.\\
Regarding the proofs, it is not hard to see that both case I and II depend only on the ratio limit result (\ref{ratio limit}),
which we believe is true for general ``critical" L\'{e}vy processes, however case III seems to be more involved.
\end{rem}

\subsection{Conditioning on large height}\label{ss:h}
In this subsection we consider the classical conditioning of large height. Assume that Assumption \ref{H2} holds and $\alpha\geq 0$,
then from Section \ref{ss:cb}
we see that $H$ is finite a.s. and has the positive continuous density $(h_x(t),t>0)$ such that
$$
h_x(t)=-x e^{-x\overline{v}_t}\frac{\partial\overline{v}_t}{\partial t}=x e^{-x\overline{v}_t}\Phi(\overline{v}_t).
$$
Now we may give the density version of the conditioning of large height.

\begin{prop}\label{h} Assume that Assumption \ref{H2} holds and $\alpha\geq 0$. Then for any $\mathcal{F}_t$-measurable bounded random variable $F$,
as $r\rightarrow\infty$,
$$
\mathbf{E}_x[F| H =r]\rightarrow \frac{1}{x}\mathbf{E}_x[ e^{\alpha t}X_t\,F]\quad \text{and}\quad
\mathbf{E}_x[F| H >r]\rightarrow \frac{1}{x}\mathbf{E}_x[ e^{\alpha t}X_t\,F].
$$
\end{prop}
\proof As in the proof of Proposition \ref{tmc}, we only need to prove the density version.
For the density version, first note that for $r>t$,
$$
\mathbf{E}_x[F\mathbf{1}\{H=r\}]=\mathbf{E}_x[F\mathbf{1}\{H-t=r-t\}].
$$
Then as in the proof of Proposition \ref{w}, by Markov property it suffices to verify that
$$
\lim_{r\rightarrow\infty}\frac{h_y(r-t)}{h_x(r)}=\frac{y}{x}e^{\alpha t},
$$
since $\mathbf{E}_x[e^{\alpha t}X_t/x]=1$.
Finally by the facts that
$\lim_{r\rightarrow\infty}\overline{v}_r=0$, $\overline{v}_r=v_t(\overline{v}_{r-t})$,
and $\frac{\partial v_t(\lambda)}{\partial \lambda}|_{\lambda=0+}=e^{-\alpha t}$, we get
$$
\lim_{r\rightarrow\infty}\frac{\partial\overline{v}_r}{\partial\overline{v}_{r-t}}=e^{-\alpha t},
$$
so
$$
\lim_{r\rightarrow\infty}\frac{h_y(r-t)}{h_x(r)}
=\frac{y}{x}\lim_{r\rightarrow\infty}\frac{\partial\overline{v}_{r-t}}{\partial\overline{v}_r}=\frac{y}{x}e^{\alpha t}.
$$
\qed

If Assumption \ref{H2} does not hold, then $H=\infty$ a.s. In this case, clearly
$\mathbf{E}_x[F| H =r]$ can not be defined and $\mathbf{E}_x[F| H >r]=\mathbf{E}_x[F]$.




\end{document}